\documentclass[11pt]{amsart}
\usepackage{amssymb}
\usepackage{amsmath}
\usepackage{amsthm}
\usepackage{eucal}

\newcommand{\eo}{\mathfrak{o}}
\newcommand{\et}{\mathfrak{t}}

\newcommand{\Ltwo}{{L}^2}

\newcommand{\Ctwo}{C_{(2)}}
\newcommand{\Uone}{U(1)}

\newcommand{\bbC}{\mathbb{C}}
\newcommand{\bbZ}{\mathbb{Z}}
\newcommand{\bbQ}{\mathbb{Q}}
\newcommand{\bbR}{\mathbb{R}}

\newcommand{\mdet}{|{\det}'|}

\DeclareMathOperator{\vertices}{Vert}
\DeclareMathOperator{\edges}{Edge}

\DeclareMathOperator{\spec}{spec}

\newtheorem{theorem}{Theorem}[section]
\newtheorem{lemma}[theorem]{Lemma}
\newtheorem{lem}[theorem]{Lemma}
\newtheorem{proposition}[theorem]{Proposition}

\newtheorem{corollary}[theorem]{Corollary}
\newtheorem{conjecture}[theorem]{Conjecture}

\theoremstyle{definition}
\newtheorem{definition}[theorem]{Definition}

\theoremstyle{remark}

\numberwithin{equation}{section}

\newcommand{\Tr}{{\text{Tr}}}

\newcommand{\im}{{\text{Im}}}

\begin{document}

%\hfill{Preliminary version: Do not circulate please}

\title[Approximating spectral invariants of Harper operators on graphs]
{Approximating Spectral invariants of Harper operators on graphs}
\author{Varghese Mathai}
\address{Department of Mathematics, MIT, Cambridge, Massachusetts 02139, USA\\
and\\
Department of Mathematics, University of Adelaide, Adelaide 5005,
Australia}
\email{vmathai@maths.adelaide.edu.au}
\author{Stuart Yates}
\address{Department of Mathematics, University of Adelaide, Adelaide 5005,
Australia}
\email{syates@maths.adelaide.edu.au}

%\date{}
\thanks{The first author acknowledges that this work was completed
in part for the Clay Mathematical Institute}

\begin{abstract}
We study Harper operators and the closely related
discrete magnetic Laplacians (DML) on a graph with a free action
of a discrete group, as defined by Sunada \cite{Sun}.
A main result in this paper is that the spectral density function
of DMLs associated to rational weight functions on graphs with a free action
of an amenable discrete group,
can be approximated by the average spectral density function
of the DMLs on a
regular exhaustion, with either Dirichlet or Neumann boundary
conditions. This then gives a criterion for the existence of
gaps in the spectrum of the DML, as well as other interesting
spectral properties of such DMLs. The technique used incorporates
some results of algebraic number theory.
\end{abstract}
\maketitle

%------------------------------------------------------------
%------------------------------------------------------------
\section*{Introduction}

Both the Harper operator and the discrete magnetic Laplacian (DML) on
the Cayley graph
of ${\mathbb Z}^2$ have been extremely well studied in mathematical physics, as
they arise as  the Hamiltonian for the
discrete model describing the quantum mechanics of free electrons in
the presence of a
magnetic field. In particular, the DML is the Hamiltonian of the
discrete model for the
integer quantum Hall effect, cf. \cite{Bel}. These operators can be
easily generalized to the Cayley graph of an arbitrary discrete group. This
and a further generalization to general graphs with a free co-compact
action of a discrete group with finite quotient, was defined by
Sunada \cite{Sun}
and studied in the context of noncommutative Bloch theory and the
fractional quantum
Hall effect in \cite{CHMM}, \cite{CHM}, \cite{MM}, \cite{MM2}.

In this paper we are primarily concerned with the approximation of the spectral
density function of the DML by examining restrictions of the DML over
a sequence of finite subgraphs.
Consider a graph $X$ which has a finite
fundamental domain under the free action of an amenable group $\Gamma$.
By F\o{}lner's characterization of amenability (see also \cite{Ad}) there
must exist a \emph{regular exhaustion} of $\Gamma$: a tower of
finite subsets $\Lambda_m\subset\Lambda_{m+1}$, $\cup_m \Lambda_m=\Gamma$
satisfying
\begin{equation}
\label{e:groupre}
\lim_{k\rightarrow\infty}\;
\frac{\#\partial_{\delta}\Lambda_{k}}{\#\Lambda_{k}}=0
\end{equation}
where
$\partial_{\delta}\Lambda_{k}=
\{\gamma\in\Gamma:d_1(\gamma,\Lambda_{k})<\delta\text{ and }
d_1(\gamma,\Gamma\setminus\Lambda_{k})<\delta\}$
is a $\delta$-neighborhood of the boundary of $\Lambda_{k}$,
and $d_1$ is the word metric on $\Gamma$ with respect to some generating set.
Given a choice $\mathcal{F}$ of fundamental domain for the $\Gamma$ action on
$X$, a corresponding sequence of subgraphs $X_m$ is constructed, with
$X_m$ being the largest subgraph of $X$ contained within
$\cup \{\gamma\mathcal{F}\,|\,\gamma\in\Lambda_m\}$. These $X_m$ satisfy a
similar property to (\ref{e:groupre}):
\begin{equation}
\label{e:graphre}
   \lim_{m\rightarrow\infty}\;
   \frac{\#{\vertices \partial_\delta X_m}}{\#{\vertices X_m}}=0
   \quad\forall\delta>0,
\end{equation}
where $\partial_\delta X_m$ refers to the subgraph which is
the $\delta$-boundary of $X_m$ in $X$, that being
the intersection of the $\delta$-neighbourhoods
of $X_m$ and $X\setminus X_m$ in the simplicial metric.

The DML restricted to the space of functions over the vertices of $X_m$
is finite dimensional, and these restricted DMLs constitute a sequence
of finite approximations to the DML itself. This leads us to the following
conjecture.

\begin{conjecture}
   Let $X$ be a graph on which there is a free group
   action by an amenable group $\Gamma$, with finite fundamental domain.
   Let $\big\{X_m\big\}^{\infty}_{m=1}$ be the regular exhaustion
   of $X$ corresponding to
   a regular exhaustion $\Lambda_m$ of $\Gamma$.  Then
   \begin{equation}
     \label{e:c0}
     \lim_{m\rightarrow\infty}\;\frac{E_m(\lambda)}{N_m}= F(\lambda),
\qquad\forall \lambda\in\bbR,
   \end{equation}
   where $N_m$ is the size of $\Lambda_m$, $E_m(\lambda)$ denotes the
   number of eigenvalues less than or equal to $\lambda$ of the DML
   restricted to
   $X_m$ with either Dirichlet or Neumann boundary conditions, and $F$
   is the spectral density function of the DML on $X$, which is defined
   using the von Neumann trace.
\end{conjecture}

We prove two approximation theorems which partially resolve this conjecture.
The first states that the equality (\ref{e:c0}) holds at all but at most
a countable set of points.

\begin{theorem}[Rough approximation theorem for spectral density functions]
   \label{t:rat}
   In the notation above, one has
   \begin{equation}
     \label{e:erat}
     \lim_{m\rightarrow\infty}\;\frac{E_m(\lambda)}{N_{m}}= F(\lambda)
   \end{equation}
   at every point of continuity of $F$.
   In particular, (\ref{e:erat}) holds for all but possibly a countable set
   of points.
\end{theorem}

Alternative proofs for this theorem exist in mathematical physics, using
for instance the `Shubin formula' (see \cite{Bel2}).
We thank J. Bellissard for bringing this reference to our attention.

Our main result is a refined approximation theorem, that holds only for DMLs
associated with rational weight functions.

\begin{theorem}[Refined approximation theorem; the rational case]
   \label{t:refat}
   $\;$ In the notation above, suppose that the DML is associated with a
   rational weight function. Then one has
   \begin{equation}
     \lim_{m\rightarrow\infty}\;\frac{E_m(\lambda)}{N_{m}}= F(\lambda),
     \qquad\forall \lambda\in\bbR.
   \end{equation}
\end{theorem}

The strategy of proof of these approximation theorems is similar to that
in \cite{DM}, but the details differ and are more involved in this case.
The key lemma for the proof of Theorem \ref{t:rat} is the combinatorial
analogue of the principle of not feeling the boundary.

Lemmas \ref{l:tr1} and \ref{l:tr2} in section 2 establish an approximation
for $F(\lambda)$ in terms of the restricted DMLs $\Delta_\sigma^{(m)}$ and
a sequence of polynomial approximations to the characteristic function
$\chi_{[0,\lambda]}$:
\begin{equation}
   \label{e:limswap}
   F(\lambda)=
   \lim_{n\rightarrow \infty}
   \lim_{m\rightarrow \infty} \frac{1}{N_m}\;
   \Tr_{{\mathbb C}}\Big(p_n\Big(\Delta_\sigma^{(m)}\Big)\Big).
\end{equation}
The proof of the refined approximation theorem \ref{t:refat} has at its
core a variant of Lemma 2.8 in \cite{Lu}. Such a result is called `Log
H\"older continuity' in the mathematical physics literature. By providing
bounds on the growth of scaled spectral densities $E_m(\lambda)/N_m$,
our variant of L\"uck's lemma allows us to effectively interchange the
limits in (\ref{e:limswap}) and provide us with our result. Showing that
the lemma is applicable is more involved in the case of the DML; we
require a positive lower bound on the modulus of the product of non-zero
eigenvalues of the restricted DMLs translated by $\lambda$. Such a bound
for algebraic $\lambda$ is established in section 3, using some
interesting facts about algebraic number fields and their rings of
integers. Proposition \ref{p:ptspec}, which states that the point
spectrum of the DML with rational weight function must be contained
within the algebraic numbers, completes the proof.

The conjecture is known to be true for all weight functions - not just
rational ones - when $X$ is the Cayley graph of $\mathbb{Z}^2$. In
this case the spectral density function of the DML can be shown to be
continuous (\cite{Sh}, \cite{DS}), and the result then follows from
Theorem \ref{t:rat}. The proof that the spectral density function is
continuous appears to utilize the special geometry of the Cayley graph
of $\mathbb{Z}^2$.  We do not know if the spectral density function is
continuous for general amenable graphs, and our proof of Theorem
\ref{t:refat} is thus quite different to that of the $\mathbb{Z}^2$
result.

In the final section of this paper we examine some consequences of
these approximation theorems. Where the weight function of the DML
is rational, Corollary \ref{c:sgc} provides a criterion for the
existence of spectral gaps which is in principle computable;
the existence of spectral gaps
for general DMLs is unknown and is a central open question in the
area.  Further in this section, we consider the Fuglede-Kadison
determinant of DML$-\lambda$, and show that it is greater than
zero for most $\lambda$.  This qualification on the spectrum of
the DML may have applications in Mathematical Physics. In joint
work with Dodziuk, Linnell and Schick \cite{DLMSY},
we have also worked out the analogs of
the results of this paper for the combinatorial Laplacian on $L^2$
cochains on covering spaces of finite $CW$ complexes.

% ----------------
% ----------------

\section{The Discrete magnetic Laplacian}

Consider a combinatorial graph $X$ on which there is a free
action by a group $\Gamma$, with finite fundamental domain. The
edge set $\edges X$ is a collection of oriented edge; each combinatorial
edge in $X$ has two corresponding oriented edges in $\edges X$, one
for each choice of orientation.

If $e$ is an oriented edge, $\overline{e}$ will denote the edge
with opposite orientation, and $\et(e)$ and $\eo(e)$ will denote the
terminus and origin respectively.

It will be convenient to regard a subset $E^+$ of $\edges X$ in which
each combinatorial edge has exactly one oriented representative;
$E^+$ corresponds to a choice of orientation for the graph.

If $g$ is a function defined on the edge set $\edges X$ with values
in some group $H$, then for compatibility with the orientation reversing
operation on the edges, we demand that $g(\overline{e})=g(e)^{-1}$.
The graph structure gives a coboundary
operator $d$ from the space of $\Ltwo$ functions on the vertices
of $X$ to the $\Ltwo$ functions on the edges, giving the chain complex
\[
0\longrightarrow\Ctwo^0(X)
\stackrel{d}{\longrightarrow}\Ctwo^1(X)\longrightarrow 0
\]
where
\[
(df)(e)=f(\et(e))-f(\eo(e)).
\]
Here $\et(e)$ denotes the terminus of the edge $e$, and $\eo(e)$ the origin.

The adjoint $d^\ast$ is given by
\[
(d^\ast g)(v) =
\sum_{\substack{e\in E^+\\\et(e)=v}}g(e)
-\sum_{\substack{e\in E^+\\\eo(e)=v}}g(e)
%&=\sum_{\substack{e\in E^+_v\\\et(e)=v}}g(e)
\]

One can then define the {\em discrete Laplacian} on $\Ctwo^0(X)$ by
\begin{gather*}
\Delta f=d^\ast df\\
(\Delta f)(v)=\mathcal{O}(v)f(v)
-\sum_{\substack{e\in E^+\\\et(e)=v}}f(\eo(e))
-\sum_{\substack{e\in E^+\\\eo(e)=v}}f(\et(e))
\end{gather*}
where $\mathcal{O}(v)$ is the valence of the vertex $v$.

Comparing this with the {\em random walk operator} on $\Ctwo^0(X)$
\[
(Rf)(v)=\sum_{\substack{e\in E^+\\\et(e)=v}}f(\eo(e))+
\sum_{\substack{e\in E^+\\\eo(e)=v}}f(\et(e)),
\]
we see that
$$
(\Delta f)(v)=\mathcal{O}(v)f(v) - (Rf)(v).
$$

The random walk operator is the basis of the generalized Harper
operator of Sunada, where the value of $f$ at each of the neighbouring
points is weighted with a complex number depending on the edge.

More precisely, consider the cochain of $U(1)$ valued functions on
the graph $X$. Two functions in $C^1(X;U(1))$ are said to be
equivalent if they belong in the same cohomology class. The
weight function $\sigma$ that we choose on the graph must
be {\em weakly $\Gamma$-invariant}, meaning that $\sigma$ is
equivalent to its left translation by any element of $\Gamma$.
Specifically, for each $\gamma\in\Gamma$ there must exist a function
$s_\gamma: \vertices X \to U(1)$ such that
\[
\sigma(\gamma e)=\sigma(e)s_\gamma(\et(e))\overline{s_\gamma(\eo(e))}.
\]
The Harper operator is then given by
\[
(H_\sigma f)(v)=
\sum_{\substack{e\in E^+_v\\ \et(e)=v}} \sigma(e)f(\eo(e))+
\sum_{\substack{e\in E^+_v\\ \eo(e)=v}} \overline{\sigma(e)}f(\et(e))
\]
and one can define the {\em discrete magnetic Laplacian} (DML) in terms
of $H_\sigma$,
\begin{equation}
(\Delta_\sigma f)(v)=\mathcal{O}(v)f(v)-(H_\sigma f)(v).
\end{equation}

The discrete Laplacian arises naturally from the cochain complex
on $X$. The DML can also be formed this way by examining a twisted
complex with coboundary operator $d_\tau$
\begin{gather*}
d_\tau:\Ctwo^0(X){\longrightarrow}\Ctwo^1(X)\\
(d_\tau f)(e)=\tau(e)f(\et(e))-\overline{\tau(e)}f(\eo(e)).
\end{gather*}
where $\tau$ is a weight function satisfying $\tau^2=\overline{\sigma}$.

While the normal discrete Laplacian commutes with the the left
$\Gamma$-translations, the DML does not. In order to discuss its
spectral properties we need to determine the appropriate von Neumann
algebra in which to examine it.

The {\em magnetic translation operators} are the maps
\[
(T_\gamma f)(x)=s_\gamma(\gamma^{-1}x)f(\gamma^{-1}x).
\]
The form a projective representation of $\Gamma$ with
\begin{equation}
\label{e:prep}
T_\gamma T_{\gamma'}=\Theta(\gamma,\gamma')T_{\gamma\gamma'} \qquad
\forall \gamma, \gamma'\in \Gamma,
\end{equation}
where $\Theta$ is a $\Uone$-valued group 2-cocycle determined by
$\sigma$ up to cohomology.

The DML commutes with these operators, and the natural context
for the examination of its spectral properties is the finite
von Neumann algebra $B(l^2(X))^{\Gamma,\sigma}$ consisting of
all the bounded operators that commute with the $T_\gamma$.
Hereafter $\Tr_{\Gamma,\sigma}$ will denote the unique trace
on this algebra, with associated von Neumann dimension
$\dim_{\Gamma,\sigma}$.

More precisely, let $A \in B(l^2(X))^{\Gamma,\sigma}$ and
$k(v_1, v_2), \; v_i \in \vertices X$ denote its kernel or matrix.
It follows that
\begin{equation}\label{onx}
s_{\gamma}(v_1) k(\gamma v_1, \gamma v_2) s_{\gamma}(v_2)^{-1} =
k(v_1, v_2)
\end{equation}
for all $\gamma\in \Gamma$ and for all $v_1, v_2 \in \vertices X$. The
von Neumann trace of $A$ is given by the expression
\begin{equation} \label{e:vNt0}
\Tr_{\Gamma,\sigma}(A) =
\sum_{v\in {\mathcal F}} k(v,v),
\end{equation}
where ${\mathcal F}$ denotes
a fundamental domain for the action of $\Gamma$ on $X$.
The expression is well defined in view of (\ref{onx}).

To express the DML
in terms of an equivariant coboundary operator, it is
necessary to define first compatible magnetic translation
operators on $\Ctwo^1(X)$.

In order to define such translations that
commute with the twisted coboundary operator, a further condition
must be imposed upon $s_\gamma$. As the sole restriction on $s_\gamma$
is that $ds_\gamma \sigma=\gamma^\ast\sigma$, one can scale it by
$k_\gamma$, constant on connected components of $X$.
The expression $s_\mu(x)s_\gamma(\mu x)\overline{s_{\gamma\mu}(x)}$
is constant on connected components of $X$ (\cite{Sun}) as is
$s_1(x)$. So setting
$k_\gamma(x)^2=\overline{s_\gamma(x)s_{\gamma^{-1}}(\gamma x)}$
and letting $s'_\gamma(x)=s_\gamma(x)k_\gamma(x)$, one has
$ds_\gamma=ds'_\gamma$ and $s'_\gamma(x)^{-1}=s'_{\gamma^{-1}}(\gamma x)$.
Hereafter it will be assumed that $s_\gamma(x)$ enjoys this
normalization property.

Recall that the weakly $\Gamma$-invariant weight function $\tau$
satisfies $\tau^2=\overline{\sigma}$. Then one can choose $t_\gamma$
with $t_\gamma^2=\overline{s_\gamma}$ and
$\tau(\gamma e)=\tau(e)t_\gamma(\et(e))\overline{t_\gamma(\eo(e))}$.
$t_\gamma$ can be extended to edges by
\[
t_\gamma(e)=\overline{t_\gamma(\et(e))t_\gamma(\eo(e))}.
\]
Extending the magnetic translation operators to $\Ctwo^1(X)$ by
\[
(T_\gamma g)(e)=t_\gamma(\gamma^{-1}e)g(\gamma^{-1}e)
\]
then allows them to commute with the twisted coboundary operator $d_\tau$,
while preserving the property \ref{e:prep}.

In order to obtain the bounds required for the proofs of
the approximation theorems in the next section, we need the following
lemmas characterizing the near-diagonality of the DML and giving
a common bound on the operator norm of the DML and its finite restrictions
$\Delta_\sigma^{(m)}$.

Let $d(v_1,v_2)$ and $d_m(v_1,v_2)$ denote the distance
between vertices $v_1$ and $v_2$ in the simplicial metric on
$X$ and $X_m$ respectively. Recall that $d(v_1,v_2)$ is the length
of the shortest sequence of edges $\{e_1,\dots,e_k\}$
with $\eo(e_1)=v_1$, $\et(e_k)=v_2$ and $\et(e_i)=\eo(e_{i+1})$ for
all $i$.

Let $D(v_1, v_2) = \left< \Delta_\sigma \delta_{v_1},
\delta_{v_2}\right>$
denote the matrix coefficient of the DML,
$\Delta_\sigma$, and
$D^{(m)}(v_1, v_2) = \left< \Delta_\sigma^{(m)}
\delta_{v_1},
\delta_{v_2}\right>$
denote the matrix coefficient of the DML,
$\Delta_\sigma^{(m)}$.

\begin{lem}
\label{l:db}
$D(v_1, v_2) = 0$ whenever $d(v_1,v_2)>1$ and
$D^{(m)}(v_1, v_2) = 0$ whenever $d_m(v_1,v_2)>1$.
There is also a positive constant $C$ independent of
$v_1, v_2$ such that
$|D(v_1, v_2)|\le C$ and $|D^{(m)}(v_1,v_2)|\le C$.
\end{lem}

Let $D^k(v_1, v_2) = \left< \Delta_\sigma^k \delta_{v_1},
\delta_{v_2}\right>$
denote the matrix coefficient of the $k$-th power of the DML,
$\Delta_\sigma^k$, and
$D^{(m)k}(v_1, v_2) = \left< \left(\Delta_\sigma^{(m)}\right)^k
\delta_{v_1},
\delta_{v_2}\right>$
denote the matrix coefficient of the $k$-th power of the DML,
$\Delta_\sigma^{(m)k}$. Then
\[
D^k(v_1, v_2) = \sum_{u_1,\ldots u_{k-1} \in \vertices X}
D(v_1, u_1)D(u_1,u_2)\ldots D(u_{k-1}, v_2)
\]
and
\[
D^{(m)k}(v_1, v_2) = \sum_{u_1,\ldots u_{k-1} \in \vertices X_m}
D^{(m)}(v_1, u_1)D^{(m)}(u_1,u_2)\ldots D^{(m)}(u_{k-1}, v_2)
\]
Then the following lemma follows easily from Lemma \ref{l:db}.

\begin{lem}
\label{l:elbnd}
Let $k$ be a positive integer. Then
$D^k(v_1, v_2) =
0$   whenever $d(v_1, v_2)>k$ and
$D^{(m)k}(v_1, v_2) = 0$ whenever $d_m(v_1,v_2)>k$.
There is also a positive constant $C$ independent of $v_1, v_2$ such that
$|D^k(v_1, v_2)|\le C^k$ and $|D^{(m)k}(v_1, v_2)|\le C^k$.
\end{lem}

Since the DML $\Delta_\sigma^k$ commutes with magnetic translations,
$\Delta_\sigma^k$ belongs to $B(l^2(X))^{\Gamma,\sigma}$ and the
von Neumann trace is given as in (\ref{e:vNt0}),
\begin {equation}
  \label{e:vNt}
  \Tr_{\Gamma,\sigma}(\Delta_\sigma^k) =
  \sum_{v\in {\mathcal F}} D^k(v,v),
\end{equation}

%it follows that
%\begin{equation}\label{inv}
%s_{\gamma}(v_1) D^k(\gamma v_1, \gamma v_2) s_{\gamma}(v_2)^{-1} =
%D^k(v_1, v_2)
%\end{equation}
%for all $\gamma\in \Gamma$ and for all $v_1, v_2 \in \vertices X$. The
%von Neumann trace of $\Delta_\sigma^k$ is given  as before by the expression
%\begin{equation} \label{vNt}
%\Tr_{\Gamma,\sigma}(\Delta_\sigma^k) =
%\sum_{v\in {\mathcal F}} D^k(v,v),
%\end{equation}
%where ${\mathcal F}$ denotes
%a fundamental domain for the action of $\Gamma$ on $X$.
%The expression is well-defined in view of (\ref{inv}).

\begin{lem}
\label{l:dmlnb}
There is a positive number $K$  such that the
operator norms of
      $\Delta_\sigma$ and of $\Delta_\sigma^{(m)}$ for all $m=1,2\ldots$
are smaller than $K^2$.
\end{lem}

\begin{proof}
The proof is similar to that in \cite{Lu}, Lemma 2.5 and uses
Lemma \ref{l:db}
together with uniform local finiteness of $X$.  More precisely we use the
fact that the valency of any vertex in $X$
is uniformly bounded,  say
$\leq b$.
\emph {A fortiori} the same is true (with the same constant $b$) for
$X_m$  for all $m$.
We now estimate the $\ell^2$ norm of
$\Delta_\sigma\kappa$ for a function $\kappa = \sum_v a_v \delta_v$
Now
$$
\Delta_\sigma \kappa =
\sum_u \left ( \sum_v D(u,v) a_v \right ) \delta_u
$$
so that
\begin{align*}
\sum_u\left|\sum_v D(u,v)a_v\right|^2
&\leq
\sum_u\left(\sum_{d(u,v)\leq 1}|D(u,v)|^2\right)
\left(\sum_{d(u,v)\leq 1} |a_v|^2\right)
\\
&\leq
C^2b\sum_u\sum_{d(u,v)\leq 1} |a_v|^2,
\end{align*}
where we have used Lemma \ref{l:db} and Cauchy-Schwartz inequality.
In the last
sum above, for every vertex $u$, $|a_u|^2$ appears at most
$b$ times.  This proves that $\|\Delta_\sigma \kappa \|^2
\leq C^2 b^2 \|\kappa \|^2$.
Identical estimate holds (with the same proof) for $\Delta_\sigma^{(m)}$
which yields the lemma if we set $K=\sqrt{C b}$.
\end{proof}

Let $\{E(\lambda):\lambda\in[0,\infty)\}$ denote the right continuous
family of spectral projections of the DML $\Delta_\sigma$.
Since $\Delta_\sigma$ commutes with magnetic translations,
so does $E(\lambda) =
\chi_{[0,\lambda]}(\Delta_\sigma)$,
for any $\lambda\in [0,\infty)$.  Let
$F:[0,\infty)\rightarrow[0,\infty)$ denote the spectral
density function, where
$$
F(\lambda)=  \Tr_{\Gamma,\sigma}(E(\lambda)).
$$
Observe that the kernel of $d$ is also given by $E(0)$, and we have
the following.

\begin{lemma}
Let $\Gamma$ be an infinite group and $X$ a connected graph. Then
$$
0 = F(0) = \dim(ker(\Delta_\sigma)).
$$
\end{lemma}

\begin{proof}
Observe that $f \in L^2(X)$ is in the kernel of the DML
$\Delta_\sigma = d_\tau^*d_\tau$
if and only if it is in the kernel of the twisted coboundary operator $d_\tau$.
That is,
\[
0=(d_\tau f)(e)=\tau(e)f(\et(e))-\overline{\tau(e)}f(\eo(e))
\]
for all edges $e \in |X|_1$. It follows that $|f(\et(e))|=|f(\eo(e))|$
for all edges $e \in |X|_1$. Since $X$ is connected, it follows that
$|f|$ is a constant function on $X$. Since $\Gamma$ is infinite,  it follows
that $f$ is identically zero.
\end{proof}

Let ${\rm spec}(\Delta_\sigma)$  denote the $L^2$-spectrum of
$\Delta_\sigma$. Then one can show that
$$
{\rm spec}(\Delta_\sigma) = \left\{\lambda\in\mathbb R : F(\lambda + \epsilon)
- F(\lambda - \epsilon) >0 \;\;\forall \epsilon>0\right\}.
$$
One can also show that $\lambda$ is an eigenvalue of
$\Delta_\sigma$ if and only if
\begin{equation}\label{point}
\liminf_{\epsilon\to 0}\left\{F(\lambda + \epsilon)
- F(\lambda - \epsilon)\right\} >0.
\end{equation}
  In fact if $\lambda$
is an eigenvalue of $\Delta_\sigma$, and $P_\lambda$ is the projection
onto the $\lambda$ eigenspace, then $$P_\lambda = \lim_{\epsilon\to 
0}\{E(\lambda + \epsilon)
- E(\lambda - \epsilon)\}.$$

%-------------------
%-------------------

\section{Proofs of the approximation theorems}

Let $E_{m}(\lambda)$ denote the number of eigenvalues $\mu$ of
$\Delta_\sigma^{(m)}$ satisfying $\mu\leq\lambda$ and which are counted with
multiplicity.

We next make the following definitions,
\begin{equation}
\label{e:fdef}
\begin{array}{lcl}
F_{m}(\lambda)& = &\displaystyle\frac{E_{m}(\lambda)}{N_{m}}\;\\[+10pt]
\overline{F}(\lambda) & = & \displaystyle\limsup_{m\rightarrow\infty}
F_{m}(\lambda) \;\\[+10pt]
\mbox{\underline{$F$}}(\lambda) & = & \displaystyle\liminf_{m\rightarrow
\infty}F_{m}(\lambda)\; \\[+10pt]
\overline{F}^{+}(\lambda) & = & \displaystyle\lim_{\delta\rightarrow +0}
\overline{F}(\lambda+\delta) \;\\[+10pt]
\mbox{\underline{$F$}}^{+}(\lambda) & = & \displaystyle\lim_{\delta
\rightarrow +0}\mbox{\underline{$F$}}(\lambda+\delta).
\end{array}
\end{equation}

The following lemma can be regarded as a combinatorial analogue of the
principle of not feeling the boundary, cf.\ \cite{DM}.

\begin{lemma}
\label{l:tr1}
$\;$Let $\Gamma$ be an amenable group and
let $p(\lambda) = \sum_{r=0}^d a_r \lambda^r$ be a polynomial.  Then,
$$
\Tr_{\Gamma,\sigma}(p(\Delta_\sigma))=
\lim_{m\rightarrow\infty}\frac{1}{N_{m}}\;
\Tr_{{\mathbb C}}\Big(p\Big(\Delta_\sigma^{(m)}\Big)\Big).
$$
\end{lemma}

\begin{proof}  First observe that if $v\in\vertices X_m$ is such that
$d(v,\partial X_m) > k$, then
$$
D^k(v, v) = \left<\Delta_\sigma^k \delta_v, \delta_v\right>
= \left<\Delta_\sigma^{(m)k} \delta_v, \delta_v\right> =
D^{(m)k}(v, v).
$$
By  (\ref{e:vNt})
$$
\Tr_{\Gamma,\sigma}(p(\Delta_\sigma))= \frac{1}{N_m} \sum_{v\in \vertices X_m}
\left< p(\Delta_\sigma)v , v\right>.
$$
Therefore we see that
\begin{multline*}
\left| \Tr_{\Gamma,\sigma}(p(\Delta_\sigma)) - \frac{1}{N_{m}}\;
\Tr_{{\mathbb C}}\Big(p\Big(\Delta_\sigma^{(m)}\Big)\Big)\right| \le
\\
\frac{1}{N_{m}} \, \sum_{r=0}^d \,  |a_r|
\sum_{\substack{v \in\vertices X_m \\ d(v, \partial X_m) \leq d}}
\, \left( |D^r(v, v)| + |D^{(m)r}(v, v)|\right).
\end{multline*}
Using Lemma \ref{l:elbnd},
we see that there is a positive constant $C$ such that
$$
\left|\Tr_{\Gamma,\sigma}(p(\Delta_\sigma)) - \frac{1}{N_{m}}\;
\Tr_{{\mathbb C}}\Big(p\Big(\Delta_\sigma^{(m)}\Big)\Big)\right| \le
2\, \sum_{r=0}^d\,\frac{\#\vertices\partial_d X_m}{N_{m}} |a_r| \, C^r.
$$
The proof of the lemma is completed by taking the limit as
$m\rightarrow\infty$, as the $X_m$ form a regular exhaustion of $X$.
\end{proof}

We next recall the following abstract lemma of L\"uck \cite{Lu}
which is proved using the Lebesgue dominated convergence theorem.

\begin{lem}
\label{l:tr2}
$\;$Let $p_{n}(\mu)$ be a sequence of polynomials
such that for
the characteristic function of the interval $[0,\lambda]$,
$\chi_{[0,\lambda]}
(\mu)$, and an appropriate real number $L$,
$$
\lim_{n\rightarrow\infty}p_{n}(\mu)=\chi_{[0,\lambda]}(\mu)\;\;\mbox{ and }
\;\;|p_{n}(\mu)|\leq L
$$
holds for each $\mu\in[0,||\Delta_\sigma||^{2}]$.  Then
$$
\lim_{n\rightarrow\infty}\Tr_{\Gamma,\sigma} (p_{n}(\Delta_\sigma))=F(\lambda).
$$\end{lem}

\noindent Lemma \ref{l:tr1} and Lemma \ref{l:tr2} give
$$
F(\lambda)  = \lim_{n\rightarrow \infty}
\Tr_{\Gamma,\sigma} (p_n(\Delta_\sigma))= \lim_{n\rightarrow \infty}
      \lim_{m\rightarrow \infty} \frac{1}{N_{m}}\;
\Tr_{{\mathbb C}}\Big(p_n\Big(\Delta_\sigma^{(m)}\Big)\Big).
$$
If one could interchange the two limits on the right-hand side above,
the proof of Theorem \ref{t:refat} would be complete.  The following lemma
is a variant
% of Lemma 2.8 due to L\"uck \cite{Lu}, and will be used in the paper
of Lemma 2.8 in \cite{Lu} and will be used in this paper
to justify
the interchange of the limits. Such a result is called `Log H\"older 
continuity'
in the mathematical physics literature, and was proved for the case 
$\Gamma = {\mathbb Z}^k$
by Craig and Simon \cite{CS}.

\begin{lem}
\label{l:ll}
$\;$Let $G:V\rightarrow V$ be a self-adjoint linear map of
the finite dimensional
Hilbert space $V$.  Let $p(t)= \det(t-G)$
be the characteristic
polynomial of $G$.  Then $p(t)$ can be factorized as
$p(t)=t^{k}q(t)$ where $q(t)$ is a polynomial with
$q(0)\neq 0$.  Let
$K$ be a real number,
$K^2\geq\max\{1,||G||\}$ and $C>0$ be a positive constant with
$|q(0)|\geq C>0$.
Let $E(\lambda)$ be the number of eigenvalues
$\eta$ of $G$, counted with multiplicity,
satisfying
$\eta\leq\lambda$.  Then for $0<\epsilon<1$,
the following
estimate is satisfied.
$$
\frac{ \{E(\epsilon)\}- \{E(0)\}}{\dim_{{\mathbb C}}V}\leq
\frac{-\log C}{(-\log\epsilon)\dim_{{\mathbb C}}V}+
\frac{\log K^{2}}{-\log\epsilon}\;.
$$
\end{lem}

Applying this lemma to the self adjoint operator $G-\lambda$,
we obtain the following estimate for $0<\epsilon<1$,
$$
\frac{ \{E(\lambda + \epsilon)\}- \{E(\lambda)\}}{\dim_{{\mathbb C}}V}\leq
\frac{-\log C(\lambda)}{(-\log\epsilon)\dim_{{\mathbb C}}V}+
\frac{\log (K^2+\lambda)}{-\log\epsilon}\;.
$$
where $C(\lambda)$ is the lower bound for $|q_\lambda(0)|$,
where $p_\lambda(t)$ is the characteristic polynomial of
$G-\lambda$, and $p_\lambda(t) = t^{k(\lambda)} q_\lambda(t)$ is its
factorization
such that $q_\lambda(0) \ne 0$.

To use this lemma, observe that $\Delta_\sigma^{(m)}$ can be regarded
as a matrix with entries in the ring of integers of the number field
generated by the weight function, whenever the weight function is
rational. Therefore we have the following

\begin{lemma}
\label{l:qmzero}
Let $\sigma$ be a rational weight function, that is $\sigma^n=1$ for some
positive integer $n$. If $\lambda$ is a non-negative algebraic number,
and $p_{m, \lambda}(t)=\det(t-(\Delta_\sigma^{(m)} - \lambda))$ is the
characteristic polynomial
of $\Delta_\sigma^{(m)}- \lambda$, and
$p_{m, \lambda}(t)=t^kq_{m, \lambda}(t)$ with
$q_{m, \lambda}(0)\neq 0$,
then there exist constants $h$ and $Q$ independent of $m$ such that
$|q_{m, \lambda}(0)|\geq (aN_m)^{-h}Q^{-haN_m}$, where $a$ is the
size of the fundamental domain $\mathcal{F}$.
\end{lemma}

\begin{proof}
See section \ref{s:ant}.
\end{proof}

In order to prove Theorem \ref{t:refat} for all $\lambda$, we will
also need the following proposition characterizing the point spectra of
the DML.

\begin{proposition}
   \label{p:ptspec}
   The point spectrum of the DML $\Delta_\sigma$ is a subset of the
   union of the spectra of the restricted DMLs $\Delta_\sigma^{(m)}$.
   In particular when $\sigma$ is a rational weight function, if there
   exists $\lambda\in\spec_{\text{\rm point}}\Delta_\sigma$, then
   $\lambda$ is algebraic.
\end{proposition}
\begin{proof}
Suppose $\lambda$ is not an eigenvalue of any of the $\Delta_\sigma^{(m)}$,
that is $\lambda\in\bbR\setminus\cup_m \spec \Delta_\sigma^{(m)}$.
Let $P_\lambda$ be the projection onto the eigenspace for $\lambda$ of
$\Delta_\sigma$. Let $X_m^o = \partial_1 X_m \cap (X\setminus X_m)$ be the
outer $1$-boundary of $X_m$ and let $P_m$ and $P_m^o$ be the projections
onto the functions with support in $X_m$ and $X_m^o$ respectively.
One has
\begin{align*}
    \Tr_{\Gamma,\sigma} P_\lambda &= \frac{1}{N_m} \Tr_{\bbC} P_mP_\lambda \\
    &\leq \frac{1}{N_m} \dim \im\, P_m P_\lambda\quad
    \text{as $\|P_m P_\lambda\|\leq 1$,}\\
    \intertext{giving}
    \Tr_{\Gamma,\sigma} P_\lambda &\leq \liminf_{m\to\infty}
    \frac{1}{N_m} \dim \im\, P_m P_\lambda.
\end{align*}
Note that by lemma \ref{l:db},
$P_m\Delta_\sigma$ differs from $\Delta_\sigma^{(m)}P_m$ only on
$\partial_1 X_m$. One can write
\begin{align*}
P_m\Delta_\sigma = P_m\Delta_\sigma P_m + B_m P_m^o \\
&= \Delta_\sigma^{(m)}P_m +B_m P_m^o
\end{align*}
where $B_m: C^0(X_m^o)\to C^0(X_m \cap \partial_1 X_m)$ encodes this
difference.

Consider $f\in\im P_m P_\lambda$, $f=P_m g$, $g$ an eigenfunction of
$\Delta_\sigma$ for $\lambda$. From $\Delta_\sigma g = \lambda g$
we have
\[
P_m\Delta_\sigma g = \lambda P_m g,
\]
and so
\begin{equation}
\label{e:bdry}
(\Delta_\sigma^{(m)} P_m + B_m P_m^o) g = \lambda P_m g.
\end{equation}
Consider two solutions $g_1$ and $g_2$ of this equation, with
$g_1|_{X_m^o}=g_2|_{X_m^o}$. Then $B_mP_m^o (g_1 - g_2)=0$ and we have
\[
\Delta_\sigma^{(m)} P_m (g_1-g_2) = \lambda P_m(g_1-g_2).
\]
As $\lambda\not\in\spec \Delta_\sigma^{(m)}$, this implies that
$P_m(g_1-g_2)=0$. So the value of $g$ on $X_m^o$ uniquely determines
$f=P_m g$, giving
\begin{align*}
\dim \im P_mP_\lambda &\leq \dim C^0(X_m^o) \\
&\leq \#\partial_1 X_m
\end{align*}
and thus
\begin{align*}
\Tr_{\Gamma,\sigma} P_\lambda
&\leq \lim_{m\to\infty} \frac{1}{N_m} \#\partial_1 X_m \\
&= 0,
\end{align*}
demonstrating that $\lambda$ is not in the point spectrum of
$\Delta_\sigma$.  As the restricted DMLs $\Delta_\sigma^{(m)}$ are all
finite operators with matrix elements belonging to the set of
algebraic numbers, the union of their spectra must in turn be a subset
of the algebraic numbers; this result therefore implies that any
$\lambda$ in the point spectrum of the DML must be algebraic.
\end{proof}

%\begin{remark}
\textbf{Remark. }
This lemma is easily generalized to any local operator $A$ in the
von Neumann algebra whose
components $A(u,v)=\langle A\delta_u,\delta_v\rangle$ are algebraic, where
we regard as local an operator $A$ for whom $A(u,v)=0$ for
all $u$ and $v$ with $d(u,v)>k$ for some fixed $k$.
\medskip
%\end{remark}

Theorem \ref{t:rat} follows immediately from part (i) of the
following theorem,
and Theorem \ref{t:refat} follows
immediately from part (ii) and the preceeding proposition \ref{p:ptspec}.

\begin{theorem}
\label{t:tmr}
Let $\Gamma$  be an amenable discrete group.

\begin{itemize}
\item[(i)]
In the notation of (\ref{e:fdef}), one has
$$
\;F(\lambda)=\overline{F}^{+}(\lambda)=
\mbox{\underline{$F$}}^{+}(\lambda).
$$
\item[(ii)] If in addition, one assumes that the weight function $\sigma$
is rational, then
$$
\;F(\lambda)=\overline{F}(\lambda)=
\mbox{\underline{$F$}}(\lambda),
$$
whenever $\lambda$ is an algebraic number. That is, under these assumptions,
one has,
$$
\;F(\lambda)=\lim_{m\to\infty}{F_m}(\lambda).
$$
\end{itemize}

\end{theorem}

\begin{proof}
Fix $\lambda\geq 0$ and define for $n\geq 1$ a continuous function
$f_{n}:{\mathbb R}\rightarrow {\mathbb R}$ by
$$
f_{n}(\mu)=\left\{\begin{array}{lcl}
1+\frac{1}{n} & \mbox{ if } & \mu\leq\lambda\\[+7pt]
1+\frac{1}{n}-n(\mu-\lambda) & \mbox{ if } &
\lambda\leq\mu\leq\lambda+\frac{1}{n} \\[+7pt]
\frac{1}{n}     & \mbox{ if } & \lambda+\frac{1}{n}\leq \mu
\end{array}\right.
$$
Then clearly $\chi_{[0,\lambda]}(\mu)<f_{n+1}(\mu)<f_{n}(\mu)$ and $f_{n}
(\mu)\rightarrow\chi_{[0,\lambda]}(\mu)$ as $n\rightarrow\infty$ for all
$\mu\in[0,\infty)$.  For each $n$, choose a polynomial $p_{n}$ such that
$\chi_{[0,\lambda]}(\mu)<p_{n}(\mu)<f_{n}(\mu)$ holds for all
$\mu\in[0,K^{2}]$.
We can always find such a polynomial by a sufficiently close approximation of
$f_{n+1}$.  Hence
$$
\chi_{[0,\lambda]}(\mu)<p_{n}(\mu)<2
$$
and
$$
\lim_{n\rightarrow\infty}p_{n}(\mu)=\chi_{[0,\lambda]}(\mu)
$$
for all $\mu\in [0,K^{2}]$.  Recall that $E_{m}(\lambda)$ denotes the number
of eigenvalues $\mu$ of $\Delta_\sigma^{(m)}$ satisfying $\mu\leq\lambda$
and counted with multiplicity.  Note that
$||\Delta_\sigma^{(m)} || \leq K^{2}$
by Lemma \ref{l:dmlnb}.
$$
\begin{array}{lcl}
\displaystyle\frac{1}{N_{m}}\;
\Tr_{{\mathbb C}}\big(p_{n}(\Delta_\sigma^{(m)})\big)
&=&\displaystyle \frac{1}{N_{m}}\sum_{\mu\in [0,K^2]}p_{n}(\mu)\\[+12pt]
\displaystyle & =&
\displaystyle\frac{ E_{m}(\lambda)}{N_{m}}+\frac{1}{N_{m}}\left\{
\sum_{\mu\in [0,\lambda ]}(p_{n}(\mu)-1)+\sum_{\mu\in (\lambda ,
\lambda + 1/n]}p_{n}(\mu)\right.\\[+12pt]
\displaystyle& &\displaystyle \hspace*{.5in}\left.+\;\sum_{\mu\in (\lambda
+ 1/n, K^2]}p_{n}(\mu)\right\}
\end{array}
$$
Hence, we see that
\begin{equation} \label{A}
F_{m}(\lambda)=\frac{ E_{m}(\lambda)}{N_{m}}
\leq\frac{1}{N_{m}}\;\Tr_{{\mathbb C}}\big(p_{n}(\Delta_\sigma^{(m)})\big).
\end{equation}
In addition,
$$
\begin{array}{lcl}
\displaystyle\frac{1}{N_{m}}\;
\Tr_{{\mathbb C}}\big(p_{n}(\Delta_\sigma^{(m)})\big)& \leq &
\displaystyle\frac{ E_{m}(\lambda)}{N_{m}}
+\;\frac{1}{N_{m}}\sup\{p_{n}(\mu)-1:\mu\in[0,\lambda]\}\;
E_{m}(\lambda) \\[+16pt]
\displaystyle &+&\displaystyle\;\frac{1}{N_{m}}\sup\{p_{n}(\mu):
\mu\in[\lambda,\lambda+1/n]\}\;
      (E_{m}(\lambda+1/n)-E_{m}(\lambda)) \\[+16pt]
\displaystyle &+&\displaystyle\;\frac{1}{N_{m}}\sup\{p_{n}(\mu):
\mu\in[\lambda+1/n,\;K^{2}]\}\;
      (E_{m}(K^{2})-E_{m}(\lambda+1/n)) \\[+16pt]
\displaystyle &\leq &\displaystyle\frac{ E_{m}(\lambda)}{N_{m}}+
\frac{ E_{m}(\lambda)}{nN_{m}}+
\frac{(1+1/n) (E_{m}(\lambda+1/n)-E_{m}(\lambda))}{N_{m}}
\\[+16pt]
\displaystyle & &\displaystyle\hspace*{.5in}+\;
\frac{(E_{m}(K^{2})-E_{m}(\lambda+1/n))}
{nN_{m}} \\[+16pt]
\displaystyle &\leq &\displaystyle
\frac{ E_{m}(\lambda+1/n)}{N_{m}}+\frac{1}{n}\;
\frac{ E_{m}(K^{2})}{N_{m}} \\[+16pt]
\displaystyle &\leq & \displaystyle F_{m}(\lambda+1/n)+\frac{a}{n}
\end{array}
$$
since $E_m(K^2)=\dim C^0 (X_m) = aN_m$ for a positive constant
$a$ independent of $m$, $a$ being the size of the fundamental domain
$\mathcal{F}$.
It follows that
\begin{equation} \label{B}
\frac{1}{N_{m}}\;\Tr_{{\mathbb
C}}\big(p_{n}(\Delta_\sigma^{(m)})\big)\leq F_{m}
(\lambda+1/n)+\frac{a}{n}.
\end{equation}
Taking the limit inferior in (\ref{B}) and the limit superior in (\ref{A}),
as $m\rightarrow\infty$, we get that
\begin{equation} \label{C}
{\overline{F}}(\lambda)\leq \Tr_{\Gamma,\sigma}\big(p_{n}(\Delta_\sigma)\big)
\leq\mbox{\underline{$F$}}(\lambda+1/n)+\frac{a}{n}.
\end{equation}
Taking the limit as $n\rightarrow\infty$ in (\ref{C}) and
using Lemma \ref{l:tr2}, we see that
$$
{\overline{F}}(\lambda)\leq F(\lambda)
\leq\mbox{\underline{$F$}}^{+}(\lambda).
$$
For all $\varepsilon>0$ we  have
$$
F(\lambda)\leq\mbox{\underline{$F$}}^{+}(\lambda)\leq\mbox{\underline{$F$}}
(\lambda+\varepsilon)\leq {\overline{F}}(\lambda+\varepsilon)
\leq F(\lambda+\varepsilon).
$$
Since $F$ is right continuous, we see that
$$
F(\lambda)={\overline{F}}^{+}(\lambda)=\mbox{\underline{$F$}}^{+}
(\lambda)
$$
proving part (i) of Theorem \ref{t:tmr}.

Next we assume that $\sigma$ be a rational weight function, that is
$\sigma^n=1$ for some positive integer $n$.
In the notation of Lemma \ref{l:qmzero}, let $\lambda$ be a non-negative
algebraic number, and
$p_{m, \lambda}(t)=\det(t-(\Delta_\sigma^{(m)} - \lambda))$
be the  characteristic polynomial
of $\Delta_\sigma^{(m)}- \lambda$, and
$p_{m, \lambda}(t)=t^k q_{m, \lambda}(t)$ with
$q_{m, \lambda}(0)\neq 0$.
Then
$|q_{m, \lambda}(0)|\geq (aN_m)^{-h}Q^{-haN_m}$,
$h$ and $Q$ being constants independent of $m$.
By Lemma \ref{l:ll} and the remarks following it,  for $\epsilon >0$ one has,
$$
\frac{F_{m}(\lambda+\epsilon)-F_{m}(\lambda)}{a}\leq
\frac{-\log ((aN_m)^{-h}Q^{-haN_m})}{(-\log\epsilon)aN_m}+
\frac{\log (K^{2} + \lambda)}
{-\log\epsilon}.
$$
That is,
\begin{equation}\label{D}
F_{m}(\lambda+\epsilon)\leq F_{m}(\lambda) -
\frac{h(\log aN_m +aN_m \log Q)}{{(\log\epsilon)} N_m}
-\frac{a\log (K^{2}+\lambda)}{\log\epsilon}.
\end{equation}
Taking limit inferior in (\ref{D}) as $m\rightarrow\infty$ yields
$$
\mbox{\underline{$F$}}(\lambda + \epsilon)\leq\mbox{\underline{$F$}}(\lambda)
-\frac{a\log Q}{\log \epsilon}
-\frac{a\log (K^{2}+\lambda)}{\log\epsilon}.
$$
Passing to the limit as $\epsilon\rightarrow +0$, we obtain
$\mbox{\underline{$F$}}(\lambda)=\mbox{\underline{$F$}}^{+}(\lambda)$.
A similar argument establishes that
${\overline{F}}(\lambda)={\overline{F}}^{+}(\lambda)$.
By part (i) of Theorem \ref{t:tmr}, one has
${\overline{F}}^{+}(\lambda)=F(\lambda)=\mbox{\underline{$F$}}^+ (\lambda)$.
Therefore we see that
$$ {\overline{F}}(\lambda) =
{\overline{F}}^{+}(\lambda)=F(\lambda)=\mbox{\underline{$F$}}^+ (\lambda)
= \mbox{\underline{$F$}}(\lambda)
$$
which
proves part (ii) of Theorem \ref{t:tmr}.
\end{proof}

% --------------------
% --------------------

\section{Determining a lower bound for $q_{m, \lambda}(0)$}
\label{s:ant}

That we can determine the lower bound in lemma \ref{l:qmzero}
depends upon the components of the matrix $\Delta_{\sigma}^{(m)}$
taking values in a ring of algebraic integers. Before providing
a proof of the lemma, some definitions and elementary results
of algebraic number theory are presented.

The weight function $\sigma$ takes values in $U(1)$. When it is
rational, that is there is some integer $n$ for which $\sigma^n=1$,
these values are $n$th roots of unity. The matrix elements in
$\Delta_{\sigma}^{(m)}$ are integer-linear combinations of the values
taken by $\sigma$, and so lie in $\bbZ[\zeta]$ where $\zeta$
is a primitive $n$th root of unity.

The natural context for examining the properties of this matrix is
the {\em algebraic number field}.

\begin{definition}
An extension $F$ of $\bbQ$ of finite degree is an
{\em algebraic number field}. The degree of the extension
$|F:\bbQ|=\dim_\bbQ F$ is called the degree of $F$. The integral
closure of $\bbZ$ in $F$ (that is, the subring of $F$ consisting
of elements $\alpha$ that satisfy $f(\alpha)=0$ for some
monic polynomial $f\in\bbZ[X]$) is termed the
{\em ring of (algebraic) integers} of $F$, and written $O_F$.
\end{definition}

For a given algebraic number $\lambda$, let $\lambda=\eta/b$ where
$\eta$ is an algebraic integer and $b\in\bbZ$. The field in which the
DML will be examined will then be $\bbQ(\zeta,\eta)$. The matrix elements
of $\Delta_\sigma^{(m)}$ are elements of $\bbZ[\zeta]$
and thus reside within the ring
of integers of this field. Let $h$ be the degree of this extension, which
will depend on the denominator $n$ of the rational weight function $\sigma$
and this algebraic integer $\eta$.

The extension $\bbQ(\zeta,\eta)$ is $\bbQ(\alpha)$ for some algebraic
$\alpha$. The minimum polynomial of $\alpha$ will be of degree $h$.
The other roots of this polynomial - the conjugates of $\alpha$ -
determine the set $\{e_i\}$ of $\bbQ$-preserving embeddings of
$\bbQ(\alpha)$ into $\bbC$.

The set of embeddings of an algebraic number field determine a norm
on that field.

\begin{definition}
Let $F$ be an algebraic number field of degree $h$, with distinct
embeddings $e_1,e_2,\ldots,e_h$ into $\bbC$. Then the norm on $F$
is defined as the product
\[
N(\alpha)=\prod_{i=1}^{h} e_i(\alpha).
\]
\end{definition}

The norm has the following important property:

\begin{theorem}
Let $N$ be the norm on an algebraic number field
$F$. Then
$N(\alpha)\in\bbQ$ for all $\alpha\in F$, and
further,
$N(\alpha)\in\bbZ$ for all $\alpha$ in the ring of integers
$O_F$.
\end{theorem}

In particular, $N(\alpha)\in\bbZ$ for all
$\alpha\in\bbZ[\zeta]$.

 From this one can arrive at the following
well-known result, presented in
\cite{Fa} and repeated
here.

\begin{lemma}
\label{l:aib}
Let $F\subset\bbC$ be an algebraic
number field of degree $h$,
with ring of integers $O_F$ and distinct
embeddings $e_1,e_2,\ldots,e_h$
into $\bbC$. Let $\alpha$ be a
non-zero element of $O_F$ satisfying
\[
|e_i(\alpha)|\leq
R\quad\forall i.
\]
Then
\[
|e_i(\alpha)|\geq R^{1-h}\quad\forall
i.
\]
\end{lemma}

\textbf{Proof of lemma \ref{l:qmzero}}.
We follow here an argument of Farber \cite{Fa}.
As discussed above, let $\{e_i\}$ be the $h$ $\bbQ$-preserving
embeddings of $\bbQ(\zeta,\eta)=\bbQ(\alpha)$ into $\bbC$, and let
$e_1$ be the identity embedding.

Now $\zeta^i=r_i(\alpha)$ for some polynomials $r_i$ of degree less
than $h$. There will be an upper bound on these co-efficients, as the
number of distinct $\zeta^i$ is finite. As there is also a bound
on the absolute values of the conjugates of $\alpha$, there is an upper
bound $R$:
\begin{align*}
|e_j(\zeta^i)| &\leq R\quad\forall i,j\\
|e_j(\eta)| &\leq R\quad\forall j
\end{align*}

Let the matrix $A_j=e_j(\Delta_\sigma^{(m)})$, applying $e_j$
component-wise. Then the matrix elements of $A_j$ satisfy the same
conditions as those of the DML in Lemma \ref{l:db} (for a different
constant $C'$) and by the same procedure by which a bound $K^2$
was determined for $\|\Delta_\sigma^{(m)}\|$, one can find $L$
such that $\|A_j\|\leq L^2$, and thus
\begin{align*}
|\Tr(A_j^i)| &\leq aN_m\|A_j\| \\
&\leq aN_m L^{2i}
\end{align*}

The $r$th symmetric polynomial of eigenvalues of $A_j$ is
therefore bounded,
\begin{align*}
|s_r(A_j)| &\leq
{{aN_m+r-1}\choose{r}}. L^{2r} \\
&\leq 2^{2aN_m-1}.L^{2r}
\qquad\text{(as $r\leq aN_m$)} \\
&< 4^{2aN_m}.L^{2aN_m} =
(4L^2)^{aN_m}.
\end{align*}
As the matrix elements of $\Delta_\sigma^{(m)}$ are algebraic
integers in $\bbQ(\alpha)$, so are the coefficients $c_i$ of
the characteristic polynomial $p_m(t)$. One has
$|e_j(c_i)|=|s_i(A_j)|<(4L^2)^{aN_m}$.

Consider $p_m(t+\lambda)=q_{m,\lambda}(t)t^k$. The coefficient
of $t^k$ in $p_m(t+\lambda)$ then is equal to $q_{m,\lambda}(0)$.
\begin{align*}
p_m(t+\lambda) &=
c_{aN_m}(t+\lambda)^{aN_m} + c_{aN_m-1}(t+\lambda)^{aN_m-1}+\cdots+c_0 \\
&= {aN_m \choose 0}c_{aN_m}\lambda^0t^{aN_m} + \\
&\quad \bigl(
{aN_m \choose 1} c_{aN_m}\lambda^1 + {aN_m-1 \choose 0}c_{aN_m-1}\lambda^0
\bigr) t^{aN_m-1} + \\
&\quad \bigl(
{aN_m\choose 2}c_{aN_m}\lambda^2+{aN_m-1\choose 1}c_{aN_m-1}\lambda^1
+{aN_m-2\choose 0}c_{aN_m-2}\lambda^0
\bigr) t^{aN_m-2} + \\
&\quad\cdots +
\\
&\quad \bigl(
{aN_m\choose aN_m}c_{aN_m}\lambda^{aN_m}+\cdots+
{0\choose 0}c_0\lambda^0
\bigr) t^0
\end{align*}
%
%
%
%
%*** START OF LARGE CORRECTED BLOCK IN SECTION 3
%
%
%
%
where $c_{aN_m}$=1. So
\[
q_{m,\lambda}(0)= \sum_{i=0}^{aN_m-k}
{i+k \choose k} c_{i+k} \lambda^i.
\]

Recall that $\lambda=\eta/b$ with $\eta$ an algebraic integer in
$\mathbb{Q}(\alpha)$. $b^{aN_m} q_{m,\lambda}(0)$ is then an
integer linear combination of the $c_i$ and powers of $\eta$ which
makes it an algebraic integer in $\mathbb{Q}(\alpha)$. As
$2^{aN_m} \geq {aN_m \choose k} \geq {aN_m-l \choose k}$, we get the
following upper bound on $|e_j(b^{aN_m} q_{m,\lambda}(0))|$:
\begin{align*}
  |e_j(b^{aN_m} q_{m,\lambda}(0))|
  &=
  \left| e_j \left( \sum_{i=0}^{aN_m-k}
      {i+k \choose k} c_{i+k} \eta^i b^{aN_m-i} \right) \right |
  \\ 
  &=
  \left| \sum_{i=0}^{aN_m-k}
    {i+k \choose k} b^{aN_m-i} e_j(\eta)^i e_j(c_{i+k}) \right |
  \\
  &\leq 
  \sum_{i=0}^{aN_m-k}
  2^{aN_m} b^{aN_m} R^{aN_m} (4L^2)^{aN_m}
  \\
  &\leq
  (8RbL^2)^{aN_m} aN_m.
\end{align*}

By lemma 3.4 then,
\[
|e_j(b^{aN_m} q_{m,\lambda}(0))|\geq
(8RbL^2)^{-haN_m}(aN_m)^{-h}\quad\forall j
\]
and in
particular
\[
|b^{aN_m} q_{m,\lambda}(0)|=|e_1(b^{aN_m} q_{m,\lambda}(0))|\geq
(8RbL^2)^{-haN_m}(aN_m)^{-h},
\]
giving
\begin{align*}
  |q_{m,\lambda}(0)| &\geq b^{-aN_m}(8RbL^2)^{-haN_m}(aN_m)^{-h} \\
  &\geq(8Rb^2L^2)^{-haN_m}(aN_m)^{-h}.
\end{align*}
%
%
%
%
%*** END OF LARGE CORRECTED BLOCK IN SECTION 3
%
%
%
%

Thus for suitable choices of constants $h$ and $Q$ independent of
$m$, one can write the inequality as
\[
|q_{m,\lambda}(0)| \geq Q^{-haN_m}(aN_m)^{-h}.
\]

In the case where $\lambda=a/b\geq 0$ is rational, one has
$|q_{m,\lambda}(0)| \geq (8B^2K^2)^{-naN_m}(aN_m)^{-n}$ where
$B$ is the maximum of $a$ and $b$, and $K^2$ is the bound
on the norm of the DML as given by lemma \ref{l:dmlnb}.

% --------------
% --------------

\section{The Fuglede-Kadison determinant of DML$-\lambda$}

In this section, we show that the Fuglede-Kadison determinant of DML$-\lambda$
is positive for most $\lambda$.
We begin though by deriving some corollaries to the main theorems
in the previous sections, using the notation there.
The following corollary should be compared to Proposition \ref{p:ptspec}.

\begin{corollary}
Let ${\rm spec}(\Delta_\sigma)$ denote the spectrum of the DML
$\Delta_\sigma$ and ${\rm spec}(\Delta_\sigma^{(m)})$ denote the
spectrum of the
DML $\Delta_\sigma^{(m)}$.
Then one has
$$
{\rm spec}(\Delta_\sigma) \subset \overline{\bigcup_{m\ge 1} {\rm
spec}(\Delta_\sigma^{(m)})}
$$
\end{corollary}

\begin{proof}
Let $\lambda_1, \lambda_2$ be points
of continuity of the spectral density function $F$ of the DML $\Delta_\sigma$
with $\lambda_1<\lambda_2$. Then by Theorem \ref{t:rat}, one has
$$
\lim_{m\to \infty} \left( F_m(\lambda_2) - F_m(\lambda_1)\right) =
F(\lambda_2) -F(\lambda_1)
$$
We also notice that
$$
F(\lambda_2) > F(\lambda_1) \Longleftrightarrow {\rm spec}(\Delta_\sigma) \cap
(\lambda_1, \lambda_2) \ne \emptyset
$$
and
$$
F_m(\lambda_2) > F_m(\lambda_1) \Longleftrightarrow {\rm
spec}(\Delta_\sigma^{(m)})
\cap (\lambda_1, \lambda_2) \ne \emptyset.
$$
This immediately implies the corollary.
\end{proof}

\begin{corollary}[Spectral gap criterion]
\label{c:sgc}
The interval $(\lambda_1, \lambda_2)$ is in a gap in the spectrum of the DML
$\Delta_\sigma$ with rational weight function $\sigma$ if and only if
$$
\lim_{m\to \infty} \left( F_m(\lambda_2) - F_m(\lambda_1)\right) = 0
$$
\end{corollary}

\begin{proof}
Notice that the interval $(\lambda_1, \lambda_2)$ is in a gap in the
spectrum of
the DML
$\Delta_\sigma$ if and only if $ F(\lambda_2) = F(\lambda_1) $. Then by
Theorem \ref{t:refat}, one has
$$
\lim_{m\to \infty} \left( F_m(\lambda_2) - F_m(\lambda_1)\right) =
F(\lambda_2) -F(\lambda_1) = 0.
$$
The converse is proved similarly to the proof above.
\end{proof}

\begin{corollary}[Spectral density estimate]
\label{c:sde}
Suppose that $\sigma$ is a rational weight function, and
    $F$ be the spectral density function of the DML $\Delta_\sigma$.
Then there are positive constants $C$ and $\delta$ such that for all
$\epsilon \in (0, \delta)$, one has for any algebraic $\mu$,
$$
F(\mu + \epsilon) -F(\mu) \le \frac{C}{-\log\epsilon}.
$$
\end{corollary}

\begin{proof}
This corollary follows immediately from Lemma \ref{l:ll}
and Lemma \ref{l:qmzero},
together with Theorem \ref{t:refat}.
\end{proof}

\subsection{Fuglede-Kadison determinants}
Here we investigate under what conditions the Fuglede-Kadison
determinant of $\Delta_\sigma-\mu$ is positive, for rational weight
function $\sigma$.
Recall that the Fuglede-Kadison determinant $\det_{\Gamma,\sigma}$
of an operator $A$
with spectral density function $G$ is defined by
\begin{equation}\label{e:FK}
\log {\det}_{\Gamma, \sigma}(A)= \int_{0^+<\lambda\leq L} \log
|\lambda| d G (\lambda)
\end{equation}
for some $L>\|A\|$,
whenever the right hand side of (\ref{e:FK}) is not
equal to $\;-\infty$, and ${\det}_{\Gamma, \sigma}(A)$
is defined to be zero if the right hand side of (\ref{e:FK})
equals  $\;-\infty$, cf. \cite{FK}. The following proposition
improves on the estimate in Corollary \ref{c:sde}.

\begin{proposition}
Let $\sigma$ be a rational weight function. Then for $\mu$
satisfying any of
\begin{enumerate}
\item
\label{p:fk1}
$\mu\not\in\spec\Delta_\sigma$,
\item
\label{p:fk2}
$\mu=0$,
\item
\label{p:fk3}
$\mu$ algebraic and $\mu\not\in\bigcup_m\spec\Delta_\sigma^{(m)}$,
\end{enumerate}
the Fuglede-Kadison determinant
$\det_{\Gamma,\sigma}(\Delta_\sigma-\mu)$ is positive.
\end{proposition}
\begin{proof}

Condition \ref{p:fk1}.
In this instance, the operator $\Delta_\sigma-\mu$ is
invertible, and positivity of the Fuglede-Kadison determinant
follows immediately.

Condition \ref{p:fk2}.
The argument for the positivity of
$\log \det_{\Gamma,\sigma} \Delta_\sigma$ exactly parallels that of
section 4 of \cite{DM}, using the lower bound on
the modified determinant
$\mdet \Delta_\sigma^{(m)}=q_{m,0}(0)$ determined in lemma
\ref{l:qmzero}.

Condition \ref{p:fk3}.
This case calls for a slightly modified version of the argument found
in section 4 of \cite{DM}, \cite{BFK}, but using the results in
sections 2 and 3 instead.

For notational convenience, let $A$ be the
DML offset by $\mu$, $A=\Delta_\sigma-\mu$, with von Neumann
spectral density
function $G(\lambda)=F(\lambda+\mu)$ where $F$ is the spectral
density function of $\Delta_\sigma$.

Let $A^{(m)}$ be the restricted DML similarly displaced by $\mu$,
$A^{(m)}=\Delta_\sigma^{(m)}-\mu$.
Recall that the \emph{normalized} spectral density functions of
$\Delta_\sigma^{(m)}$
$$
F_{m} (\lambda) = \frac 1 {N_m} E_j^{(m)} (\lambda)
$$
are right continuous, and note that the corresponding normalized
spectral density functions of $A^{(m)}$ will be given simply by
\[
G_m(\lambda)=F_m(\lambda+\mu).
\]
Observe that $G_{m}(\lambda)$ are step functions. Recall from
Lemma \ref{l:dmlnb} that there exists some $K$ such that
$\|\Delta_\sigma\|<K^2$, and $\|\Delta_\sigma^{(m)}\|<K^2$ for all $m$.
In the following, let $L$ be an algebraic number greater than $K^2+|\mu|$,
so that $\|A\|<L$ and $\|A^{(m)}\|<L$.

Denote by
$\mdet A^{(m)} $ the modified determinant of $A^{(m)}$,
that is the absolute value of the product
of all {\em nonzero} eigenvalues of $A^{(m)}$.
Choose positive $a$ and $b$ such that $a$ is less than the least
absolute value of any non-zero eigenvalue of $A^{(m)}$, and $b$ is greater
than the largest absolute value of any eigenvalue of $A^{(m)}$.
Then
\begin{equation}\label{one}
\frac 1 {N_m} \log \mdet A^{(m)} = \int_{a\leq|\lambda|\leq b }
\log |\lambda| d G_{m}
(\lambda).
\end{equation}
Integration by parts transforms this Stieltjes integral
%$\int_a^b \log \lambda d F_{m}(\lambda)$
as follows.
\begin{multline}\label{two}
\int_{a\leq|\lambda|\leq b} \log|\lambda| dG_m(\lambda)
=(\log b) \big(G_m(b)-G_m(-b)\big) \\
-(\log a) \big(G_m(a)-G_m(-a)\big)
-\int_{a\leq|\lambda|\leq b}\frac{G_m(\lambda)-G_m(0)}{\lambda}d\lambda.
\end{multline}

Integrating \ref{e:FK} by parts, one obtains
\begin{multline}\label{three}
\log {\det}_{\Gamma,\sigma}(A)=
(\log L)\big(G(L) - G(-L) \big)\\
+\lim_{\epsilon\to 0^+}
\Big\{(-\log\epsilon) \big( G(\epsilon)-G(-\epsilon)\big)
-\int_{\epsilon\leq|\lambda|\leq L}
\frac{G(\lambda)-G(0)}{\lambda}d\lambda\Big\}.
\end{multline}
Using the fact that
$\liminf_{\epsilon\to 0^+} (-\log\epsilon)
\big(G(\epsilon)-G(-\epsilon)\big) \ge 0$
one sees that
\begin{equation}\label{four}
\log {\det}_{\Gamma, \sigma} (A) \ge
(\log L) \big(G(L) - G(-L) \big) -
\int_{0^+\leq|\lambda|\leq L}
\frac {G(\lambda) - G(0)}{\lambda} d \lambda.
\end{equation}

We now complete the proof by estimating a lower bound of
$\log{\det}_{\Gamma,\sigma}(A)$ from a lower bound on
$\frac{1}{N_m}\log \mdet A^{(m)}$.
$\mdet A^{(m)}$ is the absolute value of the product of all the
non-zero eigenvalues of $A^{(m)}$, and hence for a rational
weight function $\sigma$ by lemma
\ref{l:qmzero} (and using the notation described there)
\begin{align*}
\mdet A^{(m)} &= |q_{m,\mu}(0)| \\
&\geq Q^{-haN_m}(aN_m)^{-h}
\end{align*}
and thence
\begin{equation}\label{e:logmdet}
\frac{1}{N_m}\log \mdet A^{(m)}
\geq -ha \log Q -\frac{h}{N_m} \log (aN_m)
\end{equation}
for some positive constants $Q$ and $h$ independent of $m$.

%By Lemma \ref{l:dmlnb}, there exists a positive number $K$, $1 \le K <
%\infty$, such that, for $m \ge 1$,
%$$
%||  \Delta_\sigma^{(m)} ||  \le K^2 \quad {\text{and}}\quad ||
%\Delta_\sigma ||
%\le K^2.$$

The following estimate is proved exactly as in Lemma 2.6 of \cite{DM},
so we will omit its proof here.
\begin{equation}\label{five}
\int_{0^+\leq|\lambda|\leq L} \frac{G(\lambda)-G(0)}{\lambda} d\lambda
\le
\liminf_{m \rightarrow \infty}
\int_{0^+\leq|\lambda|\leq L}
\frac{G_m(\lambda)-G_m(0)}{\lambda} d\lambda.
\end{equation}
Combining (\ref{one}) and (\ref{two}) with the inequality (\ref{e:logmdet})
we obtain
\begin{multline*}
\int_{0^+\leq|\lambda|\leq L} \frac{G_m(\lambda)-G_m(0)}{\lambda} d\lambda
\leq (\log L)\big(G_m(L)-G_m(-L)\big)\\
- \lim_{\epsilon\to 0^+} \big(\log \epsilon)(G_m(\epsilon)-G_m(-\epsilon)\big)
+ha \log Q +\frac{h}{N_m}\log (aN_m).
\end{multline*}
By the condition under regard, $\mu\not\in\spec\Delta_\sigma^{(m)}$, giving
\[
G_m(\epsilon)-G_m(-\epsilon)=F_m(\mu+\epsilon)-F_m(\mu-\epsilon)
=0\quad\text{for small $\epsilon$},
\]
and so one has
\begin{equation}\label{six}
\int_{0^+\leq|\lambda|\leq L} \frac{G_m(\lambda)-G_m(0)}{\lambda} d\lambda
\leq (\log L)\big(G_m(L)-G_m(-L)\big)
+ha \log Q +\frac{h}{N_m}\log (aN_m).
\end{equation}

 From (\ref{four}), (\ref{five}) and (\ref{six}), we conclude that
\begin{multline}\label{seven}
\log {\det}_{\Gamma, \sigma} \Delta_\sigma \ge (\log L)
\big( G(L)-G(-L) \big)
- ha \log Q\\
- \liminf_{m \rightarrow \infty} \Big\{
(\log L) \big(G_m(L)-G_m(-L) \big) + \frac{h}{N_m}\log (aN_m) \Big\}
\end{multline}
Now $h/N_m \log (aN_m)\to 0$ and by part (ii) of Theorem \ref{t:tmr}
\begin{align*}
G(L) = F(L+\mu) &= \lim_{m\to\infty} F_m(L+\mu)\\
&=\lim_{m\to\infty} G_m(L)
\end{align*}
(and similarly for $-L$), hence
\[
\log{\det}_{\Gamma,\sigma} A\geq -h \log Q,
\]
that is,
\[
{\det}_{\Gamma,\sigma} (\Delta_\sigma-\mu)>0.
\]

\end{proof}

\end{document}